\documentclass[]{article}
\begin{document}
\newtheorem{proposition}{Proposition}[section]
\newtheorem{definition}{Definition}[section]
\newtheorem{lemma}{Lemma}[section]

\title{\bf Algebraic Regularity over Quaternions\\and\\Regular Four-Manifolds}
\author{Keqin Liu\\Department of Mathematics\\The University of British Columbia\\Vancouver, BC\\
Canada, V6T 1Z2}
\date{November, 2015}
\maketitle

\begin{abstract} Based on a new generalization of Cauchy-Riemann system presented in this paper, we introduce a class of quaternion-valued functions of a quaternionic variable, which are called algebraic regular functions. The set of algebraic regular functions is not only a real associative algebra, but also respect the composition of functions. Using  algebraic regular functions as transition maps, we introduce a class of four-manifolds called  the regular four-manifolds.
\end{abstract}

\medskip
\section{Introduction}

The glorious history of complex analysis and the beauty of the quaternions have stimulated people to explore the analysis over the quaternions since the beginning of the last century. It is unquestionable that a successful theory of the analysis over the quaternions has to come from a right way of introducing the differentiability for quaternion-valued functions of a quaternionic variable. A lot of attempts have been made to get a right counterpart of the complex  differentiability in the context of the quaternions. For example, in 1935, R. Fueter \cite{F} introduced a class of quaternion-valued functions called Fueter regular functions by using a generalization of Cauchy-Riemann system. The theory of Fueter regular functions has been developed so well that many fundamental results about complex analytic functions have been rewritten in the quaternions \cite{S}. In 2006, G. Gentili and D.C. Struppa
\cite{GS} introduced slice regular functions of a quaternionic variable by means of the algebraic property that the quaternion algebra is a union of a family of complex planes. The development of the theory of the slice regular functions is so rapid that a monograph \cite{GSS} about the  slice regular functions was published in 2013. It is no doubt that both Fueter regular functions  and  slice regular functions are very interesting quaternion-valued functions of a quaternionic variable, but each of the two classes of functions of a quaternionic variable does not respect the composition of functions. This odd algebraic feature of Fueter regular functions  and  slice regular functions  prevent them from playing the role of transition maps in the exploration of the manifolds over the quaternions. The purpose of this paper is to present a way of fixing this odd algebraic feature. In section 2, based on a new generalization of Cauchy-Riemann system, we introduce a class of quaternion-valued functions of a quaternionic variable, which are called algebraic regular functions in this paper. The set of algebraic regular functions is not only a real associative algebra, but also respect the composition of functions. In section 3, using  algebraic regular functions as transition maps, we initiate the study of a class of four-manifolds, which are called regular four-manifolds by us.

\medskip
\section{Algebraic Regularity}

Let $\mathcal{H}=\mathcal{R}e_1\oplus \mathcal{R}e_2\oplus\mathcal{R}e_3\oplus\mathcal{R}e_4$ be the quaternions discovered by W. R. Hamilton in 1843, where $\mathcal{R}$ is the real number field, $e_1$ is the identity of the real division associative algebra $\mathcal{H}$, and the multiplication among the remaining three elements in the $\mathcal{R}$-basis
$\{e_1,\, e_2,\, e_3, \,e_4\}$ is defined by
$$ e_2^2=e_3^2=e_4^2=-e_1, \qquad e_ie_j=-e_je_i=(-1)^{i+j+1}e_{9-i-j},$$
where $2\le i<j\le 4$.

\medskip
We now introduce another way of defining  differentiability in the context of the quaternions.

\medskip
\begin{definition}\label{def2.1} Let $U$ be an open subset of $\mathcal{H}$. We say that a quaternion-valued function $f: U\to \mathcal{H}$ is {\bf algebraic regular} at
$a=\displaystyle\sum_{i=1}^4 a_ie_i\in U$ if $f$ has two properties given below.
\begin{description}
\item[(i)] There exist two $C^1$ real-valued functions $f_0: \mathcal{R}^4\to \mathcal{R}$ and $f_1: \mathcal{R}^4\to \mathcal{R}$ such that
\begin{equation}\label{eq1}
f(q)=f_1(q_1,\, q_2,\, q_3,\, q_4)e_1+
\displaystyle\sum_{i=1}^4 q_if_0(q_1,\, q_2,\, q_3,\, q_4)e_i
\end{equation}
for all $q=\displaystyle\sum_{i=1}^4 q_ie_i\in U$.
\item[(ii)] The following equations hold at $(a_1,\, a_2,\, a_3, \,a_4)\in \mathcal{R}^4$:
\begin{equation}\label{eq2}
\frac{\partial f_1}{\partial q_1}=f_0+q_2\frac{\partial f_0}{\partial q_2}+
q_3\frac{\partial f_0}{\partial q_3}+q_4\frac{\partial f_0}{\partial q_4},
\end{equation}
\begin{equation}\label{eq3}
\frac{\partial f_1}{\partial q_i}=-q_i\frac{\partial f_0}{\partial q_1}, \qquad\qquad
q_i\frac{\partial f_0}{\partial q_j}=q_j\frac{\partial f_0}{\partial q_i},
\end{equation}
where $2\le i,\, j\le 4$ and $i\ne j$.

We say that  $f: U\to \mathcal{H}$ is an {\bf algebraic regular function} on $U$ if $f$ is algebraic regular at every point of $U$.
\end{description}
\end{definition}

\medskip
Algebraic regular functions appear naturally among quaternion-valued functions of a quaternionic variable. For example, the power function $f(q)=q^n$ is an algebraic regular function on $\mathcal{H}$ for a nonnegative constant integer $n$, and $f(q)=q^n$ is an algebraic  regular function on $\mathcal{H}\setminus\{0\}$ for a negative constant integer
$n$. Another natural example of algebraic regular functions is the exponential function
$e^q$ of a  a quaternionic variable $q$, where $e^q$ is defined by
$$ e^q=\left(e^{q_1}\cos\,\sqrt{q_2^2+q_3^2+q_4^2}\right)e_1+
\sum_{i=2}^4q_i
\left(\frac{e^{q_1}\,\sin\,\sqrt{q_2^2+q_3^2+q_4^2}}{\sqrt{q_2^2+q_3^2+q_4^2}}\right)e_i$$
for $q=\displaystyle\sum_{i=1}^4 q_ie_i$ with $q_2^2+q_3^2+q_4^2\ne 0$. It is easy to check that the exponential function: $q\to e^q$ is an algebraic  regular function on
$\mathcal{H}\setminus \{q_1e_1\,|\, q_1\in \mathcal{R}\}$.

\medskip
We call the system consisting of (\ref{eq2}) and (\ref{eq3}) the {\bf generalized Cauchy-Riemann system}. The next proposition gives an equivalent description of the generalized
Cauchy-Riemann system in terms of Fueter operators.

\medskip
\begin{proposition}\label{pr2.1} Let  $U$ be an open subset of $\mathcal{H}$. If a function $f: U\to \mathcal{H}$ satisfies (\ref{eq1}), then the following are equivalent:
\begin{description}
\item[(i)] $f$ is algebraic regular at $a=\displaystyle\sum_{i=1}^4 a_ie_i\in U$;
\item[(ii)] $\mathcal{D}_{\ell}\big(f(q)\big)|_{q=a}=
\mathcal{D}_{r}\big(f(q)\big)|_{q=a}=-2f_0(a_1,\, a_2,\, a_3,\, a_4)e_1$,
where $\mathcal{D}_{\ell}=\displaystyle\sum_{i=1}^4 e_i\,\frac{\partial }{\partial q_i}$ and
 $\mathcal{D}_{r}=\displaystyle\sum_{i=1}^4 \frac{\partial }{\partial q_i}\,e_i$ are the {\bf left Fueter operator} and the {\bf right Fueter operator}, respectively.
\end{description}
\end{proposition}
\hfill\raisebox{1mm}{\framebox[2mm]{}}

\medskip
The following proposition implies that the set of algebraic regular functions on an open subset of the quaternions is an associative algebra over the real number field.

\medskip
\begin{proposition}\label{pr2.2} Let  $U$ be an open subset of $\mathcal{H}$, and let $c$ be a real number. Suppose that $f: U\to \mathcal{H}$ and $g: V\to \mathcal{H}$ are algebraic regular at $a\in U$. Then
\begin{description}
\item[(i)] The real scalar product $cf$, the sum $f+g$ and the pointwise product $fg$ are algebraic regular at $a$.
\item[(ii)] If $f(a)\ne0$, then the reciprocal $\displaystyle\frac{1}{f}$ is algebraic regular at $a$.
\end{description}
\end{proposition}
\hfill\raisebox{1mm}{\framebox[2mm]{}}

\medskip
The most significant property of algebraic regular functions is the following
\medskip
\begin{proposition}\label{pr2.3} Let  $U$ and $V$ be  open subsets of $\mathcal{H}$, and let $f: U\to \mathcal{H}$ and $g: V\to \mathcal{H}$ be two quaternion-valued functions of a quaternionic variable, where $f(U)\subseteq V$. If $f$ is algebraic regular at $a\in U$ and
$g$ is algebraic regular at $f(a)\in V$, then the composition function $g\circ f$ is algebraic regular at $a$.
\end{proposition}
\hfill\raisebox{1mm}{\framebox[2mm]{}}

\medskip
\section{Regular Four-Manifolds}

By identifying $\mathcal{H}$ with $\mathcal{R}^4$ by means of the map:
$\displaystyle\sum_{i=1}^4 q_ie_i\to (q_1,\,q_2,\,q_3,\,q_4)$, we introduce a class of four-manifolds in the following

\medskip
\begin{definition}\label{def3.1} Let $M$ be a smooth four-manifold with a $C^{\infty}$-atlas
$\{(U_{\alpha},\, \phi_{\alpha})\,|\, \alpha\in \mathcal{S}\}$, where
$\phi_{\alpha}: U_{\alpha}\to\mathcal{H}$ is a homeomorphism from the open subset $U_{\alpha}$ of $M$ onto an open subset $\phi_{\alpha}( U_{\alpha})$ of $\mathcal{H}$ for all $\alpha\in \mathcal{S}$. The  smooth four-manifold $M$ is called a {\bf regular four-manifold} if the transition map
$$
\phi_{\alpha}\circ\phi_{\beta}^{-1}: \phi_{\beta}(U_{\alpha}\cap U_{\beta})\to
\phi_{\alpha}(U_{\alpha}\cap U_{\beta})
$$
is an algebraic regular function on the open subset $\phi_{\beta}(U_{\alpha}\cap U_{\beta})$ of
$\mathcal{H}$ whenever for $\alpha$,$\beta\in \mathcal{S}$ and $U_{\alpha}\cap U_{\beta}\ne \emptyset$.
\end{definition}

\medskip
We finish this paper with two examples of regular four-manifolds.

\medskip
{\bf Example 1} The quaternions $\mathcal{H}$ is a regular four-manifold with a $C^{\infty}$-atlas which consists of a unique chart $(id_{\mathcal{H}},\, \mathcal{H})$, where $id_{\mathcal{H}}$ is the identity function on the quaternions $\mathcal{H}$.

\medskip
{\bf Example 2} On the set $\mathcal{H}\times \mathcal{H}\setminus \{(0,\,0)\}$, there is an equivalence relation $\sim$ defined by
$$
(x_1,\, x_2)\sim (y_1,\, y_2)\,\iff
\,\mbox{$y_1=x_1q$ and $y_2=x_2q$ for some $q\in \mathcal{H}\setminus \{(0\}$,}
$$
where $(x_1,\, x_2)$, $(y_1,\, y_2)\in \mathcal{H}\times \mathcal{H}\setminus \{(0,\,0)\}$. Let
$[(x_1,\, x_2)]$ be the equivalence class containing $(x_1,\, x_2)$, and let
$$
\mathcal{H}\mathcal{P}:=\{\, [(x_1,\, x_2)]\,|\,
(x_1,\, x_2)\in \mathcal{H}\times \mathcal{H}\setminus \{(0,\,0)\}\}.
$$
The quotient topology induced from $\mathcal{H}\times \mathcal{H}\setminus \{(0,\,0)\}$ with respect to the
equivalence relation $\sim$ makes $\mathcal{H}\mathcal{P}$ a topological space. For $i=1,\, 2$, let
$$
U_i:=\{\, [(x_1,\, x_2)]\,|\,
(x_1,\, x_2)\in \mathcal{H}\times \mathcal{H}\setminus \{(0,\,0)\}\quad\mbox{and}\quad x_i\ne 0\},
$$
and define $\phi_i: U_i\to \mathcal{H}$ by
$$ \phi_1([(x_1,\, x_2)]): =x_2x_1^{-1}\quad\mbox{for $[(x_1,\, x_2)]\in U_1$}$$
and
$$ \phi_2([(x_1,\, x_2)]): =x_1x_2^{-1}\quad\mbox{for $[(x_1,\, x_2)]\in U_2$}.$$
Then  $\mathcal{H}\mathcal{P}$ is a regular four-manifold with the
$C^{\infty}$-atlas $\{(U_i,\, \phi_i\,|\, i=1,\, 2\}$, which is called the {\bf right projective regular four-manifold} by us.


\bigskip
\begin{thebibliography}{9}
\bibitem{F} R. Fueter, \textsl{Die Funktionentheorie der Differentialeichungen $\Delta u=0$ and
$\Delta\Delta u=0$ mit vier reellen Variablen}, Comment. Math. Helv. 7 (1935), 307-330
\bibitem{GS} G. Gentili and D.C. Struppa, \textsl{A new approach to Cullen-regular functions of a quaternionic variable}, C. R. Math. Acda. Sci. Paris 342 (10) (2006), 741-744
\bibitem{GSS} G. Gentili, C. Stoppato and D.C. Struppa, \textsl{Regular Functions of Quaternionic Variables}, Springer Monographs in Mathematics, Springer-Verlag Berlin Heidelbeig, 2013
\bibitem{S} A. Sudbery, \textsl{Quaternionic analysis}, Math. Proc. Camb. Phil. Soc. 85 (1979), 199-225
\end{thebibliography}
\end{document}